\documentclass[12 pt]{article}
\usepackage{amscd, amssymb,amsmath,amsthm}
\newtheorem{lm}{Lemma}

\newtheorem{Theorem}{Theorem}

\newcommand{\ext}{\operatorname{Ext}}
\newcommand{\Hom}{\operatorname{Hom}}
\newcommand{\sheafext}{\mathcal{E}xt}

\newcommand{\T}{\mathbf{T}}
\newcommand{\A}{\mathbf{A}}
\newcommand{\C}{\mathbf{C}}
\newcommand{\proj}{\mathbf{P}}
\newcommand{\Z}{\mathbb{Z}}
\newcommand{\R}{\mathbb{R}}
\newcommand{\N}{\mathcal{N}}
\newcommand{\cO}{\mathcal{O}}
\newcommand{\cI}{\mathcal{I}}
\newcommand{\cT}{\mathcal{T}}
\newcommand{\bV}{\mathsf{V}}
\newcommand{\bE}{\mathsf{E}}
\newcommand{\bw}{\mathsf{w}}
\newcommand{\bv}{\mathsf{v}}

\newcommand{\bF}{\mathsf{F}}
\newcommand{\bZ}{\mathsf{Z}}

\DeclareMathOperator{\tr}{tr}

\newcommand{\rarr}{\rightarrow}
\newcommand{\oh}{{\mathcal{O}}}

\newcommand{\eqq}{\stackrel{\sim}{=}}

\newcommand{\bpf}{\noindent {\em Proof.} }
\newcommand{\epf}{\qed \vspace{+10pt}}

\begin{document}

\title{Gromov-Witten theory and Donaldson-Thomas theory, I}
\author{D.~Maulik, N.~Nekrasov, A.~Okounkov, and R.~Pandharipande}
\date{1 May 2004}

\maketitle

\vspace*{-8.5cm}
\hfill{\scriptsize{ITEP-TH-61/03}} 

\hfill{\scriptsize{\ IHES/M/03/67}} 
\vspace*{8.5cm}

\section{Introduction}

\subsection{Overview}

Let $X$ be a nonsingular, projective, Calabi-Yau 3-fold.
Gromov-Witten theory concerns counts of maps of curves to $X$.
The counts are defined in terms of a canonical $0$-dimensional 
virtual fundamental class on the moduli space of maps. 
The discrete invariants of a map are the genus $g$  of the
domain and the  
 degree $\beta\in H_2(X, \Z)$ of the image.
For every $g$ and $\beta$, the Gromov-Witten invariant
is the virtual number of genus $g$, degree $\beta$ maps. 
We sum the contributions of all genera with weight $u^{2g-2}$, 
where $u$ is a parameter. 

The Gromov-Witten invariants have long been expected to 
be expressible in terms of
appropriate curve counts in the target $X$. 
A curve in $X$ corresponds to an ideal sheaf on $X$.
The discrete invariants of the ideal sheaf are the holomorphic 
Euler characteristic $\chi$
and the fundamental class $\beta\in H_2(X,\Z)$ of the associated
curve.
Donaldson and Thomas have constructed a canonical $0$-dimensional 
virtual fundamental class on the moduli space of ideal sheaves on $X$.
For every $\chi$ and $\beta$, the Donaldson-Thomas invariant is the virtual number
of the corresponding ideal sheaves. We  sum the contributions  
over $\chi$ with weight $q^\chi$, 
where $q$ is a parameter. 

We present here a precise mathematical conjecture relating the Gromov-Witten
and Donaldson-Thomas theories of $X$. Our conjecture is 
motivated by the description of Gromov-Witten theory via crystals in \cite{ORV}.
A connection between Gromov-Witten theory and
integration over the moduli space of ideal sheaves is strongly suggested
there. A related physical conjecture is 
formulated in \cite{harvard}.

\vspace{+10pt}
\noindent {{\bf Conjecture.}}
{\em  The change of variables,
$
e^{iu}=-q
$,
equates the Gromov-Witten and Donaldson-Thomas theories of $X$.}
\vspace{+10pt}

The moduli of maps and sheaves have been
related previously by the Gopakumar-Vafa conjecture
equating Gromov-Witten invariants
to BPS state counts determined by
the {\em cohomology} of the moduli of sheaves \cite{GoVafa1,GoVafa2}.
The Gopakumar-Vafa conjecture has been verified in several cases
and has been a significant source of motivation. However,
there have been
difficulties on the mathematical side in selecting an appropriate
cohomology theory for the singular moduli of sheaves which arise, see \cite{HST}.

Donaldson-Thomas theory concerns {\em integration} over the
moduli of sheaves. The subject was defined, along with a construction
of the virtual class, by Donaldson and Thomas in \cite{DonThom,Thom}
with motivation from several sources, see \cite{AOS,BKS,witten}.
As the Donaldson-Thomas invariant is similar to the Euler
characteristic of the moduli of sheaves, a philosophical
connection between Gromov-Witten
invariants and the cohomology of the moduli of sheaves is implicit
in our work.
However, the Donaldson-Thomas invariant is {\em not}
the Euler characteristic.

As evidence for our conjecture, we present a proof in the toric local
Calabi-Yau
case via the virtual localization formula for
Donaldson-Thomas theory. The proof depends upon partial 
evaluations of
the topological vertex on the Gromov-Witten side.

\subsection{General 3-folds}
We believe the Gromov-Witten/Donaldson-Thomas 
correspondence holds for all 3-folds.
Donaldson-Thomas theory has a natural supply of
observables constructed from the Chern classes of universal
sheaves. These Chern classes should
correspond to insertions in Gromov-Witten
theory, see \cite{LMN}. 
The degree $0$ case, where no insertions are required,
is discussed in Section \ref{degz} below. For primary fields,
a complete GW/DT correspondence for
all 3-folds is conjectured in \cite{MNOP2}. For descendent fields,
the correspondence is precisely formulated for the
descendents of a point in \cite{MNOP2}.

We conjecture
the \emph{equivariant vertex} \cite{N2},
discussed briefly in Section \ref{seqv} below,
has the same relation to general cubic Hodge integrals as the
topological vertex does to Calabi-Yau Hodge integrals \cite{DF}.
A closely related issue is the precise formulation of the
GW/DT correspondence for {\em all} descendent fields.
We will pursue the topic in a future paper.

\subsection{GW theory}
Gromov-Witten theory is defined via integration over the moduli
space of stable maps.
Let $X$ be a nonsingular, projective, Calabi-Yau 3-fold.
Let $\overline{M}_g(X,\beta)$ denote the moduli space of
stable maps from connected genus $g$ curves to $X$ representing the
class $\beta\in H_2(X, \Z)$, and
let
$$N_{g,\beta} = \int_{[\overline{M}_g(X,\beta)]^{vir}} 1,$$
denote the corresponding
Gromov-Witten invariant. Foundational aspects of the theory
are treated, for example, in \cite{Beh, BehFan, LiTian}.

Let $\bF'_{GW}(X;u,v)$ denote the reduced Gromov-Witten potential of $X$,
$$
\bF'_{GW}(X;u,v)= \sum_{\beta\neq 0} \sum_{g\geq  0} N_{g,\beta}
\ u^{2g-2} v^\beta,
$$
omitting the constant maps.
The reduced partition function,
$$
\bZ'_{GW}(X;u,v) = \exp \bF'_{GW}(X;u,v)\,,
$$
generates disconnected Gromov-Witten invariants of $X$
excluding constant contributions.

Let
$\bZ'_{GW}(X;u)_\beta$
denote the reduced partition function of degree $\beta$
invariants,
$$
\bZ'_{GW}(X;u,v) = 1+
\sum_{\beta\neq 0} \bZ'_{GW}(X;u)_\beta \ v^\beta.$$

\subsection{{DT theory}}
Donaldson-Thomas theory is defined via integration over
the moduli space of ideal sheaves.
Let $X$ be a nonsingular, projective, Calabi-Yau 3-fold.
An {\em ideal sheaf} is a torsion-free sheaf of rank 1 with
trivial determinant. Each ideal sheaf ${\mathcal I}$
injects into its double dual,
$$0 \rarr {\mathcal I} \rarr {\mathcal I} ^{\vee\vee}.$$
As ${\mathcal I}^{\vee\vee}$ is reflexive of rank 1 with trivial
determinant,
$${\mathcal I}^{\vee\vee}\eqq \oh_X,$$
see \cite{OSS}.
Each ideal sheaf $\cI$
determines a subscheme $Y\subset X$,
$$0 \rarr {\mathcal I} \rarr \oh_X \rarr \oh_{Y} \rarr 0.$$
The maximal dimensional components of $Y$
(weighted by their intrinsic multiplicities) determine
an element,
$$[Y] \in H_*(X,\Z).$$
Let $I_n(X,\beta)$ denote the moduli space of ideal sheaves
${\mathcal I}$ satisfying
$$\chi(\oh_{Y}) = n,$$
and
$$[Y] = \beta\in H_2(X,\Z).$$
Here, $\chi$ denotes the holomorphic Euler characteristic.
The moduli space $I_n(X,\beta)$ is isomorphic to the
Hilbert scheme of curves of $X$ \cite{MP}.

The Donaldson-Thomas invariant is defined
via integration against the dimension 0 virtual class,
$$\tilde{N}_{n,\beta} = \int_{[I_n(X,\beta)]^{vir}} 1.$$
Foundational aspects of the theory
are treated in \cite{MP,Thom}.

Let $\bZ_{DT}(X;q,v)$ be the partition function of the Donaldson-Thomas
theory of $X$,
$$
\bZ_{DT}(X;q,v)=
\sum_{\beta\in H_2(X,\Z)}
\sum_{n\in \Z} \tilde{N}_{n,\beta} \ q^n \, v^\beta\,.
$$
An elementary verification shows, for fixed $\beta$, the invariant
$\tilde{N}_{n,\beta}$ vanishes for sufficiently negative $n$ since
the corresponding moduli spaces of ideal sheaves are empty.

The degree 0 moduli space $I_n(X,0)$ is
isomorphic to the Hilbert scheme of $n$ points on $X$.
The degree 0 partition function,
$$
\bZ_{DT}(X;q)_0=
\sum_{n\geq 0} \tilde{N}_{n,\beta} \ q^n\,,
$$
plays a special role in the theory.
The McMahon function,
$$
M(q) = \prod_{n\geq 1} \frac{1}{(1-q^n)^n}\,,
$$
is the generating series for $3$-dimensional partitions,
see \cite{stanley}.

\vspace{+10pt}
\noindent {{\bf Conjecture 1.}}  The degree 0 partition function is
determined by
$$\bZ_{DT}(X;q)_0= M(-q)^{\chi(X)},$$
where $\chi(X)$ is the topological Euler characteristic.
\vspace{+10pt}

The reduced partition function $\bZ'_{DT}(X;q,v)$ is defined by
quotienting by the degree 0 function,
$$
\bZ'_{DT}(X;q,v) = \bZ_{DT}(X;q,v) \big/ \bZ_{DT}(X;q)_0
\,.
$$
Let $\bZ'_{DT}(X;q)_\beta$ denote the reduced partition function of degree
$\beta\neq 0$ invariants,
$$
\bZ'_{DT}(X;q,v)= 1+
\sum_{\beta\neq 0}  \bZ'_{DT}(X;q)_\beta \, v^\beta.$$

\vspace{+10pt}
\noindent {{\bf Conjecture 2.}}
The reduced series
$\bZ'_{DT}(X;q)_\beta$ is a rational function of $q$
symmetric under the transformation $q\mapsto 1/q$.
\vspace{+10pt}

We now state our main conjecture relating the Gromov-Witten theory and the
Donaldson-Thomas
theory of a Calabi-Yau 3-fold $X$.

\vspace{+10pt}
\noindent {{\bf Conjecture 3.}} The change of variables $e^{iu}=-q$
equates the reduced partition functions:
$$
\bZ'_{GW}(X;u,v) = \bZ'_{DT}(X;-e^{iu},v) \,.
$$
\vspace{+5pt}

The change of variables in Conjecture 3 is well-defined by Conjecture 2.
Gromov-Witten and Donaldson-Thomas theory
 may be viewed as expansions of a single partition
functions at different points. Conjecture 3 can be checked order by order
in $u$ and $q$ only if an effective bound on
the degree of the rational function in Conjecture 2 is known.

\subsection{Integrality}
The Gopakumar-Vafa conjecture
for Calabi-Yau 3-folds
predicts the following form for  the reduced
Gromov-Witten partition function via the change of variables $e^{iu}=-q$,
$$
\bZ'_{GW}(X;u)_\beta = q^r \frac{f(q)}{\prod_{i=1}^k (1-(-q)^{s_i})^2}, \ \ \
f\in\Z[q],\ \ \, r\in \Z, \ \  s_i>0\,.
$$
In particular, by the Gopakumar-Vafa conjecture, $\bZ'_{GW}(X;u)_\beta$
defines a series in $q$ with {\em integer} coefficients.

Conjecture 3 identifies the $q$ series with the reduced
partition function of the Donaldson-Thomas theory of $X$.
Integrality of the Donaldson-Thomas invariants holds by
construction (as no orbifolds occur).
We may refine Conjecture 2 above to fit the form of
the Gopakumar-Vafa conjecture.

\subsection{Gauge/string dualities}

In spirit, our conjecture is similar to a 
gauge/string duality with Donaldson-Thomas theory on the gauge side and 
 Gromov-Witten 
theory on the string side. 
While at present we are not aware of a purely 
gauge theoretic interpretation of  
Donaldson-Thomas theory, there are various
indications such an interpretation should exist. Most
importantly, the equivariant vertex measure
which appears in the equivariant localization formula for
Donaldson-Thomas theory, see Section \ref{seqv}, is 
identical to the equivariant vertex 
in noncommutative Yang-Mills theory \cite{N1}. 
We plan to investigate the issue 
further.

The interplay between gauge fields and strings is one of the
central themes in modern theoretical and mathematical
physics \cite{Polyakov}. In particular, the conjectural Chern-Simons/string
duality of Gopakumar and Vafa \cite{GoV} was a source of
many insights into the Gromov-Witten theory of Calabi-Yau
3-folds. As a culmination of these developments, the
\emph{topological vertex} was introduced in
\cite{topver}. The topological vertex is a certain explicit
function of three partitions $\lambda$, $\mu$, $\nu$, and the genus expansion
parameter $u$, which is an elementary building block
for constructing the Gromov-Witten invariants of arbitrary
local toric Calabi-Yau 3-folds. The gauge/string duality seems to hold in a
broader context, see \cite{LMN}, \cite{icmtalk} for evidence in the
Fano case.

In \cite{ORV}, the topological
vertex was interpreted as counting 3-dimensional
partitions $\pi$ with asymptotics $\lambda,\mu,\nu$
along the coordinate axes. The variable $q=e^{iu}$
couples to the volume of the partition $\pi$ in
the enumeration. 
The global data obtained by gluing such 3-dimensional
partitions according to the gluing rules of the topological
vertex  was observed in \cite{harvard,ORV} to naturally corresponds to torus invariant
ideal sheaves in the target 3-fold $X$. 
The main mathematical result of our paper is the identification
of the topological vertex expansion with the
equivariant  localization
formula for the Donaldson-Thomas theory of the local Calabi-Yau
geometry.

The GW/DT correspondence is conjectured
to hold for {\em all} Calabi-Yau 3-folds. While several motivations
for the correspondence came from local Calabi-Yau geometry, new
methods of attack will be required to study the full GW/DT correspondence.

A relation between Gromov-Witten theory and gauge theory on
the same space $X$ has been observed
previously in four (real) dimensions in the context of
Seiberg-Witten invariants \cite{taubes}. There, a deformation of the Seiberg-Witten
equations by a 2-form yields solutions concentrated near the zero locus of the 2-form,
an embedded curve. 

We expect, in our case, the sheaf-theoretic
description of curves will be identified with 
a deformed version of solutions to some
gauge theory problem. An outcome 
 should be a natural method of deriving the
Donaldson-Thomas measure. The gauge theory in question is a deformation of
the twisted maximally supersymmetric Yang-Mills theory compactified on our 3-fold
$X$ \cite{N1}. The theory, discussed in \cite{BKS}, has BPS solutions and generalized instantons. The
expansion of the super-Yang-Mills action about these solutions gives rise to a
quadratic form with bosonic and fermionic determinants which should furnish the
required measure \cite{ihiggs1, ihiggs2}. 

In case $X = {\C}^3$, the deformation in question is
the passage to the noncommutative ${\R}^6$, see \cite{noncom1, noncom2}.
Ordinary gauge theories have typically non-compact moduli spaces of BPS
solutions. It is customary in mathematics to compactify these spaces by replacing 
holomorphic bundles by coherent sheaves. The physical consequences of such a replacement 
are usually quite interesting
and lead to many insights \cite{avatars1, avatars2,avatars3}.
Sometimes the ``compactified''
space is non-empty while the original space is empty. Our problem corresponds to $U(1)$ gauge 
fields which do not support
nontrivial instantons, while the compactified moduli space of instantons is non-empty and 
coincides with the
Hilbert scheme of curves of given topology
on $X$.

\subsection{Acknowledgments}

We thank J. Bryan,
A. Iqbal, M. Kontsevich, Y. Soibelman, R. Thomas, and C. Vafa for related discussions.

D.~M. was partially supported by a Princeton Centennial graduate fellowship.
A.~O.\ was partially supported by
DMS-0096246 and fellowships from the Sloan and Packard foundations.
R.~P.\ was partially supported by DMS-0071473
and fellowships from the Sloan and Packard foundations.

\section{Degree 0}
\label{degz}

\subsection{GW theory}
Let $X$ be a nonsingular, projective 3-fold (not necessarily Calabi-Yau).
The
degree 0 potential  $\bF_{GW}(X;u)_0$ may be separated as:
$$
\bF_{GW}(X;u)_0 =
\bF^0_{X,0} + \bF^1_{X,0} + \sum_{g\geq 2} \bF^g_{X,0}\,.
$$
The genus 0 and 1  contributions in degree 0 are not constants,
the variables of the classical cohomology appear explicitly.
Formulas can be found, for example, in \cite{icmtalk}.

We will be concerned here with the higher genus terms. For $g\geq 2$,
a virtual class calculation yields,
$$
\bF^g_{X,0} = (-1)^g\frac{u^{2g-2}}
{2}
 \int_X \Big( c_3(X) -  c_1(X)c_2(X) \Big)         \ \cdot
\int_{\overline{M}_g} \lambda_{g-1}^3,$$
where $c_i$ and $\lambda_i$ denote the Chern classes of the tangent bundle $T_X$ and
and the Hodge bundle ${\mathbb E}_g$ respectively.
Define the degree 0 partition function of Gromov-Witten theory by
$$
\bZ_{GW}(X;u)_0 = \exp\left(\sum_{g\geq 2} \bF^g_{X,0}\right)\,.
$$
The Hodge integrals which arise have been computed in \cite{FabPan},
\begin{equation}
  \label{BB}
  \int_{\overline{M}_g}\lambda_{g-1}^3=\frac{|B_{2g}|}{2g}
\frac{|B_{2g-2}|}{2g-2}\frac1{(2g-2)!},
\end{equation}
where $B_{2g}$ and $B_{2g-2}$ are Bernoulli numbers.

Using the Euler-Maclaurin formula, 
the asymptotic relation,
\begin{equation}
  \label{GWMM}
  \bZ_{GW}(X;u)_0 \sim M(e^{iu})^{\frac{1}{2}
\int_X  c_3(X) -  c_1(X)c_2(X) } \,,
\end{equation}
may be derived from \eqref{BB}.
The precise meaning of \eqref{GWMM} is the following:
the logarithms of both sides have identical $o(1)$-tails
in their $u\to 0$ asymptotic expansion.

\subsection{DT theory}
We now turn to the degree 0 partition function for the Donaldson-Thomas
theory of $X$. The first issue is the construction of the virtual class
in Donaldson-Thomas theory.

In \cite{Thom}, the Donaldson-Thomas theory of $X$
is defined only in the Calabi-Yau and Fano cases. 
In fact, the
arguments of \cite{Thom} use only the existence of an anticanonical section
on $X$.

\begin{lm} \label{p}
Let $X$ be a nonsingular, projective 3-fold
with $$H^0(X,\wedge^3 T_X) \neq 0,$$ then $I_n(X,\beta)$
carries a canonical perfect obstruction theory.
\end{lm}

Under the hypotheses of Lemma \ref{p}, the Donaldson-Thomas theory
of $X$ is constructed for {\em higher} rank sheaves as well as the rank 1 case
of ideal sheaves.
The connection, if any, between 
Gromov-Witten theory and the higher rank
Donaldson-Thomas theories is not clear to us.

The technical
condition required in \cite{Thom} for the construction of the perfect
obstruction theory and the virtual class $[I_n(X,\beta)]^{vir}$
is the vanishing of traceless
$\ext^3_0({\mathcal I}, {\mathcal I})$
for all $[{\mathcal I}]\in I_n(X,\beta)$.
In fact, $\ext^3_0({\mathcal I}, {\mathcal I})$ vanishes for 
every ideal sheaf on a nonsingular, projective 3-fold $X$ \cite{MP}.
Hence, Donaldson-Thomas theory is well-defined in rank 1 for {\em all} $X$.

For simplicity, let us assume the vanishing of the higher cohomology of
the structure sheaf,
\begin{equation}
\label{jjjf}
H^i(X,\oh_X)=0,
\end{equation}
for $i\geq 1$. Then, $\ext_0({\mathcal I}, {\mathcal I})$ equals $\ext({\mathcal I}, {\mathcal I})$.

\begin{lm} \label{pp} Let $X$ be a nonsingular, projective, 3-fold satisfying \eqref{jjjf}.
Then,
$$\ext^3({\mathcal I},{\mathcal I})=0,$$
for all $[{\mathcal I}] \in I_n(X,0)$.
\end{lm}

\bpf
By Serre duality for $\text{Ext}$, 
$$\ext^3({\mathcal I},{\mathcal I})=\ext^0({\mathcal I},{\mathcal I}\otimes K_X)^\vee,$$
where $K_X$ denotes the canonical bundle. We must therefore prove
$$\text{Hom}({\mathcal I},{\mathcal I}\otimes K_X)=0.$$ 
Let $U\subset X$ be the complement of the support of $Y$.
Since ${\mathcal I}$ restricts to $\oh_U$ on $U$,
$$\text{Hom}({\mathcal I}|_U,{\mathcal I}|_U\otimes K_U)=\Gamma(U,K_U)=H^0(X,K_X).$$ 
The last equality is obtained from the extension of sections since
$Y$ has at most 1-dimensional support.
Since ${\mathcal I}$ is torsion-free, the restriction,
$$\text{Hom}({\mathcal I},{\mathcal I}\otimes K_X) \rarr
\text{Hom}({\mathcal I|_U},{\mathcal I|_U}\otimes K_U),$$
is injective.
Since $h^0(X,K_X)=h^3(X,\oh_X)$, the Lemma is proven.
\epf

 The proof of the vanshing of
$\ext^3_0({\mathcal I},{\mathcal I})$
in the presence of higher cohomology of the
structure sheaf is similar, see \cite{MP}.

The virtual dimension of $I_n(X,0)$ is 0 for general 3-folds $X$.
A simple calculation from the definitions yields the following result.

\begin{lm} \label{ppp} $\tilde{N}_{1,0} = -\int_X c_3(X)- c_1(X)c_2(X).$
\end{lm}

\bpf 
The moduli space $I_{1}(X,0)$
is the nonsingular 3-fold $X$.
The tangent bundle is $\ext_0^1({\mathcal I, \mathcal I})$, and
the obstruction bundle is $\ext_0^2({\mathcal I}, {\mathcal I})$. Using
Serre duality and the local-to-global spectral sequence for $\ext$, we find
the obstruction bundle is isomorphic to $(T_X\otimes K_X)^\vee$. Then,
$$\tilde{N}_{1,0} = -\int_X c_3( T_X \otimes K_X ) =
-\int_X c_3(X)-c_1(X)c_2(X), $$
completing the proof.
\epf

The degree 0 Gromov-Witten and Donaldson-Thomas
theories are already related by  Lemma \ref{ppp}. However, we make a
stronger connection generalizing Conjecture 1.

\vspace{+10pt}
\noindent {{\bf Conjecture $\mathbf 1'$.}}
The degree 0 Donaldson-Thomas partition function for a 3-fold $X$ is
determined by:
$$
\bZ_{DT}(X;q)_0= M(-q)^{\int_X c_3(T_X \otimes K_X)}\,.
$$
\vspace{+10pt}

We will present a proof of Conjecture $1'$ in case $X$ is a nonsingular
toric 3-fold in \cite{MNOP2}.

The series $M(q)$ arises naturally in the computation of the
Euler characteristic of the Hilbert scheme of points
of a 3-fold \cite{cheah}. It would be interesting to find a
direct connection between the degree 0 
Donaldson-Thomas invariants and the Euler characteristics of $I_n(X,0)$ 
in the Calabi-Yau case.

\section{Local Calabi-Yau geometry} \label{ldef}

\subsection{GW theory}
Let $S$ be a nonsingular, projective, toric,
Fano surface with canonical bundle $K_S$.
The Gromov-Witten theory of the
local Calabi-Yau geometry of $S$ is defined via an excess
integral.
Denote the universal curve and universal map over the
moduli space of stable maps to $S$ by:
$$\pi: U \rarr \overline{M}_g(S,\beta),$$
$$\mu: U \rarr S.$$
Then,
$$N_{g,\beta}= \int_{[\overline{M}_g(S,\beta)]^{vir}} e( R^1\pi_*\mu^* K_S),$$
for $0 \neq \beta \in H_2(S,\Z)$.
The reduced partition function $\bZ'_{GW}(X;u,v)$ is defined in terms
of the local invariants $N_{g,\beta}$ as before.

\subsection{DT theory}
Let $X$ be the projective bundle $\proj(K_S\oplus \oh_S)$
over the surface $S$. The Donaldson-Thomas theory of $X$ is well-defined
in every rank by the following observation.

\begin{lm} $X$ has an anticanonical section.
\end{lm}

\bpf
Consider the fibration $\pi: X \rarr S$.
We have,
$$\wedge^3 T_X = T_\pi \otimes \pi^*( \wedge^2 T_S),$$
where $T_\pi$ is the $\pi$-vertical tangent line.

Let $V$ denote the vector bundle $K_S \oplus \oh_S$ on $S$.
The $\pi$-relative Euler sequence is:
$$0 \rarr \oh_X \rarr \pi^*(V) \otimes \oh_{\proj(V)}(1) \rarr T_\pi \rarr 0.$$
Hence,
$$T_\pi = \wedge^2 \pi^*(V) \otimes \oh_{\proj(V)}(2).$$
we conclude,
$$\wedge^3 T_X =  \wedge^2 \pi^*(V) \otimes \oh_{\proj(V)}(2) \otimes \pi^*( \wedge^2 T_S)=
\oh_{\proj(V)}(2).$$
However, since
$$H^0(X,\oh_{\proj(V)}(2)) = H^0(S, Sym^2 V^*) \neq 0,$$
the Lemma is proven.
\epf

For classes $\beta\in H_2(S,\Z)$, we define the
reduced partition function for the Donaldson-Thomas theory of the
Calabi-Yau geometry of $S$ by
$$
\bZ'_{DT}(S;q)_\beta=\bZ'_{DT}(X;q)_\beta \,.
$$
While $X$ is not Calabi-Yau, the Donaldson-Thomas theory of
$X$ is still well-defined by Lemma \ref{p}.

We will prove Conjectures 2 and 3 are true for the local Calabi-Yau
geometry of toric Fano surfaces by virtual localization.

\subsection{Local curves}
The constructions above also define the local Calabi-Yau theory of
the curve $\proj^1$ with normal bundle $\oh(-1) \oplus \oh(-1)$.
The proof of Conjectures 2 and 3 for local surfaces given in Section \ref{lcal} below is valid
for the local $\proj^1$ case.

The Gromov-Witten theory of a local Calabi-Yau curve of arbitrary genus 
has been defined in \cite{BryPan}. We believe the
GW/DT correspondence holds for these geometries as well \cite{BP1}.

\section{Localization in Donaldson-Thomas theory}
\label{lcal}
\subsection{Toric geometry}
\label{cech}
Let $X$ be a nonsingular, projective, toric 3-fold.
Let $\mathbf{T}$ be the 3-dimensional complex torus
acting on $X$.
Let $\Delta(X)$ denote the
Newton polyhedron of $X$ determined by a polarization.
The polyhedron $\Delta(X)$ is the image of $X$ under
the moment map.

The vertices of $\Delta(X)$ correspond to fixed points
$X^\T=\{X_\alpha\}$
of the $\T$-action. For each $X_\alpha$, there is a canonical,
$\T$-invariant, affine open chart,
$$U_\alpha\cong \A^3,$$ centered
at $X_\alpha$. We may choose coordinates $t_i$ on
$\T$ and coordinates $x_i$ on
$U_\alpha$ for which
the $\T$-action on $U_\alpha$
is
determined by 
\begin{equation}
  \label{st_act}
  (t_1,t_2,t_3)\cdot x_i = t_i x_i  \,.
\end{equation}
In these coordinates,
the tangent representation $X_\alpha$
 has character $$t^{-1}_1+ t_2^{-1}+ t_3^{-1}.$$
We will use the covering $\{ U_\alpha\}$ of $X$ to compute
Cech cohomology.

The $\T$-invariant lines of $X$ correspond to the edges of $\Delta(X)$.
More precisely, if
$$C_{\alpha\beta}\subset  X$$
is a $\T$-invariant line incident to the fixed points
$X_\alpha$ and $X_\beta$, then
$C_{\alpha\beta}$ corresponds to an edge of $\Delta(X)$
joining the vertices $X_\alpha$ and $X_\beta$.

The geometry of $\Delta(X)$ near the edge is determined by the normal
bundle $\N_{C_{\alpha\beta}/X}$. If
$$
\N_{C_{\alpha\beta}/X} = \cO(m_{\alpha\beta}) \oplus \cO(m'_{\alpha\beta})
$$
then the transition functions between the charts $U_\alpha$ and
$U_\beta$ can be taken in the form
\begin{equation}\label{trans}
  (x_1,x_2,x_3) \mapsto (x_1^{-1}, x_2 \,x_1^{-m_{\alpha\beta}} ,
x_3 \,x_1^{-m'_{\alpha\beta}} ) \,.
\end{equation}
The curve $C_{\alpha\beta}$ is defined in these coordinates by $x_2=x_3=0$.

\subsection{Moduli of ideal sheaves}
The ${\mathbf T}$-action on $X$ canonically lifts to the moduli space
of ideal sheaves $I_n(X,\beta)$. The perfect obstruction theory
constructed by Thomas is canonically ${\mathbf T}$-equivariant \cite{MP,Thom}.
The virtual localization formula reduces integration against
$[I_n(X,\beta)]^{vir}$ to
a sum fixed point contributions \cite{GrabPan}.

The first step is to determine the ${\mathbf T}$-fixed points
of $I_n(X,\beta)$. If $$[\mathcal I]\in I_n(X,\beta)$$ is ${\mathbf T}$-fixed,
then the associated subscheme
$Y\subset X$ must be preserved
 by the torus action. Hence, $Y$
 must be supported on the ${\mathbf T}$-fixed points $X_\alpha$ and
the ${\mathbf T}$-invariant lines connecting
them.

Since $\mathcal{I}$ is $\mathbf T$-fixed on each
open set,
${\mathcal I}$  must
be defined on $U_\alpha$ by a monomial ideal,
$$
I_\alpha = \cI\big |_{U_\alpha} \subset \C[x_1,x_2,x_3] \,,
$$
and may
also be viewed as a $3$-dimensional partition $\pi_\alpha$
$$
\pi_\alpha = \left\{(k_1,k_2,k_3),\, \prod_{1}^3 x_i^{k_i} \notin
I_\alpha\right\}\subset \Z_{\ge 0}^3 \,.
$$
The  ideals $I_\alpha$ are $1$-dimensional.
The corresponding partitions $\pi_\alpha$ may be infinite
in the direction of the coordinate axes.

The asymptotics of $\pi_\alpha$ in the coordinate directions
are described by three ordinary $2$-dimensional partitions.
In particular, in the direction of the $\T$-invariant curve $C_{\alpha\beta}$,
we have the partition $\lambda_{\alpha\beta}$ with the
following diagram:
\begin{align*}
 \lambda_{\alpha\beta}&=\left\{(k_2,k_3),  \forall k_1  \, \,
x_1^{k_1} x_2^{k_2} x_3^{k_3}  \notin I_{\alpha}\right\} \\
&=\left\{(k_2,k_3),  \,
x_2^{k_2} x_3^{k_3}  \notin I_{\alpha\beta}\right\}\,,
\end{align*}
where
$$
I_{\alpha\beta} = \cI\big |_{U_\alpha\cap U_\beta}
\subset \C[x_1^{\pm 1},x_2,x_3] \,.
$$

In summary,  a $\T$-fixed ideal sheaf $\cI$ can be
described in terms of the following data:
\begin{itemize}
\item[(i)] a $2$-dimensional partition $\lambda_{\alpha\beta}$
assigned to each edge of $\Delta(X)$,
\item[(ii)] a $3$-dimensional partition $\pi_\alpha$ assigned to each
vertex of $\Delta(X)$ such that the asymptotics of
$\pi_\alpha$ in the three coordinate directions is given
by the partitions $\lambda_{\alpha\beta}$ assigned to the corresponding
edges.
\end{itemize}

\subsection{Melting crystal interpretation}\label{smelt}

The partition data $\{\pi_\alpha,\lambda_{\alpha\beta}\}$ corresponding
to a $\T$-fixed ideal sheaf $\cI$
can be given a melting crystal interpretation \cite{ORV}.
Consider the weights of the $\T$-action on
$$
H^0(X,\oh_{Y}(d))\,.
$$
For large $d$, the corresponding
points of $\Z^3$ can be described as follows.

Scale the Newton polyhedron $\Delta(X)$ by a
factor of $d$. Near each corner of $d\Delta(X)$, the
intersection $\Z^3 \cap d \Delta(X)$ looks like
a standard $\Z_{\ge 0}^3$, so we can place the corresponding
partition $\pi_\alpha$ there. Since $d$ is large
 and since, by construction $\pi_\alpha$ and $\pi_\beta$
agree along the edge joining them, a global combinatorial
object emerges.

One can imagine the points of $\Z^3 \cap d \Delta(X)$
are atoms in a crystal and, as
the crystal is melting or dissolving, some of the atoms
near the corners and along the edges are missing.
These missing atoms are described by the partitions $\pi_\alpha$ and $\lambda_{\alpha\beta}$.
They are precisely the weights of the  $\T$-action on
$H^0(X,\oh_{Y}(d))$.

\subsection{Degree and Euler characteristic}\label{sdegch}

Let $[\cI]\in I_n(X,\beta)$ be a $\T$-fixed ideal sheaf on $X$
described by the partition
 data $\{\pi_\alpha,\lambda_{\alpha\alpha'}\}$.
We see
$$
\beta = \sum_{\alpha,\alpha'} |\lambda_{\alpha\alpha'}| \, [C_{\alpha\alpha'}] \,,
$$
where $|\lambda|$ denotes the size of a partition $\lambda$, the number
of squares in the diagram.

For $3$-dimensional partitions $\pi$ one can similarly define
their size $|\pi|$ by the number of cubes in their diagram.
Since the partitions $\pi_\alpha$ may be infinite along
the coordinate axes, the number $|\pi_\alpha|$ so defined
will often be infinite. We define
the renormalized volume $|\pi_\alpha|$ as follows. Let
$\lambda_{\alpha\beta_i}$, $i=1,2,3$, be the asymptotics of
$\pi_\alpha$. We set
$$
|\pi_\alpha| = \#
\left\{\pi_\alpha \cap [0,\dots,N]^3 \right\}-
(N+1) \sum_1^3 |\lambda_{\alpha\beta_i}| \,, \quad N\gg 0 \,.
$$
The renormalized volume is
independent of the
cut-off $N$ as long as $N$ is
sufficiently large. The number $|\pi_\alpha|$
so defined may be negative.

Given $m,m'\in\Z$ and a partition $\lambda$, we
define
$$
f_{m,m'}(\lambda) = \sum_{(i,j)\in \lambda} \left[ m(i-1) + m'(j-1)+1\right] \,.
$$
Each edge of $\Delta(X)$ is assigned a pair of integers $(m_{\alpha\beta},m'_{\alpha\beta})$
from the normal bundle of the associated $\T$-invariant line and a partition
$\lambda_{\alpha\beta}$ from  the $\T$-fixed ideal sheaf $\cI$.
By definition, we set
\begin{equation}
  \label{fab}
  f(\alpha,\beta) = f_{m_{\alpha\beta},m'_{\alpha\beta}} (\lambda_{\alpha\beta}) \,.
\end{equation}

\begin{lm}\label{fchi}
$\chi(\oh_Y) = \sum_\alpha |\pi_\alpha| + \sum_{\alpha,\beta}
f(\alpha,\beta) \,.$
\end{lm}

\bpf
The result is an elementary calculation in toric geometry.
For example, a computation of $\chi(\oh_Y)$ using the
Cech cover defined by $\{U_\alpha\}$ immediately
yields the result.
\epf

\subsection{The $\T$-fixed obstruction theory}
The moduli space $I_n(X,\beta)$ carries a $\T$-equivariant perfect
obstruction theory,
$$E_0 \rarr E_1,$$
see \cite{MP,Thom}.
Assume the virtual dimension of $I_n(X,\beta)$ is 0.
The virtual localization formula \cite{GrabPan} may be stated as follows,
$$\int_{[I_n(X,\beta)]^{vir}} 1 =
\sum_{[\cI] \in I_n(X,\beta)^\T} \int_{[S(\cI)]^{vir}} \frac{e(E^m_1)}{e(E^m_0)}.$$
Here, $S(\cI)$ denotes the $\T$-fixed subscheme of $I_n(X,\beta)$
supported at the point $[\cI]$,
and $E_0^m$, $E_1^m$ denote the nonzero $\T$-weight spaces.
The virtual class, $[S(\cI)]^{vir}$, is determined
by the $\T$-fixed obstruction theory.

We first prove $S(\cI)$ is
the reduced point $[\cI]$. It suffices to prove the Zariski tangent
space to $I_n(X,\beta)$ at $[\cI]$ contains no trivial
subrepresentations.
Since $X$ is toric, all the higher cohomologies of $\oh_X$ vanish,
$$H^i(X,\oh_X)=0,$$
for $i\geq 0$.
Hence, the traceless condition is satisfied, and the Zariski  tangent
space of $I_n(X,\beta)$ at  $[{\mathcal I}]$ is
$\ext^1({\mathcal I}, {\mathcal I})$.

\begin{lm}
\label{ppww}
Let $[{\mathcal I}] \in I_n(X,\beta)$ be a ${\mathbf T}$-fixed
point. The ${\mathbf T}$-representation,
$$\ext^1(\cI,\cI),$$
contains no trivial subrepresentations.
\end{lm}

\bpf From the ${\mathbf T}$-equivariant ideal sheaf sequence,
\begin{equation}
\label{idd}
0 \rarr {\mathcal I} \rarr \oh_X \rarr \oh_{Y} \rarr 0,
\end{equation}
we obtain a sequence of ${\mathbf T}$-representations,
$$\rarr \ext^0({\mathcal I}, \oh_{Y}) \rarr \ext^1({\mathcal I}, {\mathcal I})
\rarr \ext^1({\mathcal I}, \oh_X) \rarr.$$
The left term, $\ext^0({\mathcal I}, \oh_{Y})$, does not
contain trivial representations by Lemma \ref{www} below.

We will prove the right term, $\ext^1({\mathcal I}, \oh_X)$, also
does not contain trivial representations.
By Serre duality, it
suffices to study the representation,
$$\ext^2(\oh_X,{\mathcal I}\otimes K_X)= H^2(X,{\mathcal I}\otimes K_X).$$
The long exact sequence in cohomology obtained from (\ref{idd}) by tensoring
with $K_X$ and the vanishings,
$$H^1(X,K_X)=H^2(X,K_X)=0,$$
together yield a ${\mathbf T}$-equivariant isomorphism,
$$H^1(X, \oh_{Y}\otimes K_X) \stackrel{\sim}{\rarr} H^2(X, {\mathcal I}\otimes
K_X).$$
The first Cech cohomology of $\oh_{Y}\otimes K_X$ is computed
via the representations
$$H^0(U_{\alpha\beta}, \oh_{Y}\otimes K_X),$$
where $U_{\alpha\beta}=U_\alpha \cap U_\beta.$
Here we use the Cech cover defined in Section \ref{cech}.
An elementary argument shows these representations contain
no trivial subrepresentations.
\epf

\begin{lm} \label{www}
$\ext^k({\mathcal I}, \oh_{Y})$ contains no
trivial subrepresentations.
\end{lm}
\bpf
By the local-to-global spectral sequence, it suffices to prove
$$H^i(\sheafext^j(\cI, \oh_{Y}))$$
contains no trivial subrepresentations for all $i$ and $j$.
By a Cech cohomology calculation, it then suffices to prove
$$H^0(U_\alpha, \sheafext^j(\cI_\alpha, \oh_{Y_\alpha})), \
H^0(U_{\alpha\beta}, \sheafext^j(\cI_{\alpha\beta}, \oh_{Y_{\alpha\beta}}))
$$
contain no trivial subrepresentations.
Triple intersections need not be considered since $\oh_{Y_{\alpha\beta\gamma}}$
vanishes.

We will study $Ext^j(\cI_\alpha,\oh_{Y_\alpha})$ on $U_\alpha$
via the $\T$-equivariant Taylor resolution of the
monomial ideal $\cI_\alpha$. The argument for $Ext^j(\cI_{\alpha\beta}, \oh_{Y_{\alpha\beta}})$
on $U_{\alpha\beta}$ is identical.

Let $\cI_\alpha$ be generated
by the monomials $m_1,\ldots, m_s$. For each subset
$$T \subset \{1, \dots, s\},$$ let
$$m_{T} = x^{r(T)}=\text{least common multiple of }\{m_{i} | i \in T\}. $$
For $ 1 \leq t \leq s$, let $F_{t}$ be the free $\Gamma(U_\alpha)$-module 
with basis $e_{T}$ indexed by
subsets $T \subset \{1,\dots,s\}$
of size $t$.

  A differential $d: F_{t} \rightarrow F_{t-1}$ is defined as follows.  Given a subset $T$ of
size $t$, let $T = \{i_{1}, \dots, i_{t}\}$ where $i_{1} < \dots < i_{t}$.  Let
$$d(e_{T}) = \sum_{T' = T \backslash \{i_{k}\}} (-1)^{k} x^{r_{T} - r_{T'}}.$$  The Taylor resolution,
$$ 0 \rightarrow F_{s} 
\rightarrow \dots \rightarrow F_{2} \rightarrow F_{1} \rightarrow \cI_\alpha \rightarrow 0,$$
is exact \cite{Taylor}.  
Moreover, the resolution is equivariant with $\T$-weight $r(T)$ on the
generator $e(T)$.

The weights of the generators of $F_t$ are weights of monomials in $\cI_\alpha$.
However, the weights of the $\T$-representation
$\oh_{Y_\alpha}$ are precisely not equal to weights of
monomials in $\cI_\alpha$. Hence,
$Hom(F_t,\oh_{Y_\alpha})$ contains no trivial subrepresentations.
We then conclude $Ext^j(\cI_\alpha, \oh_Y)$ contains no trivial
subrepresentations by computing via the Taylor resolution of $\cI_\alpha$.
\epf

The obstruction space at $[\cI]\in I_n(X,\beta)$ of the perfect obstruction
theory is $\ext^2(\cI,\cI)$.
The following Lemma implies the $\T$-fixed obstruction theory at $[\cI]$
is trivial.

\begin{lm}
$\ext^2(\cI,\cI)$
contains no trivial subrepresentations.
\end{lm}

\bpf From the ${\mathbf T}$-equivariant ideal sheaf sequence,
we obtain,
$$\rarr \ext^1({\mathcal I}, \oh_{Y}) \rarr \ext^2({\mathcal I}, {\mathcal I})
\rarr \ext^2({\mathcal I}, \oh_X) \rarr.$$
The left term, $\ext^1({\mathcal I}, \oh_{Y})$, does not
contain trivial representations by Lemma \ref{www} above.

We will prove the right term, $\ext^2({\mathcal I}, \oh_X)$, also
does not contain trivial representations.
By Serre duality, it
suffices to study the representation,
$$\ext^1(\oh_X,{\mathcal I}\otimes K_X)= H^1(X,{\mathcal I}\otimes K_X).$$
The long exact sequence in cohomology obtained from  by tensoring the ideal
sheaf sequence
with $K_X$ and the vanishings,
$$H^0(X,K_X)=H^1(X,K_X)=0,$$
together yield a ${\mathbf T}$-equivariant isomorphism,
$$H^0(X, \oh_{Y}\otimes K_X) \stackrel{\sim}{\rarr} H^1(X, {\mathcal I}\otimes
K_X).$$
The space of global sections of $\oh_{Y}\otimes K_X$ is computed
via the representations
$$H^0(U_\alpha, \oh_{Y}\otimes K_X).$$
As before, an elementary argument shows these representations contain
no trivial subrepresentations.
\epf

The virtual localization formula may then be written as
$$\int_{[I_n(X,\beta)]^{vir}} 1 =
\sum_{[\cI] \in I_n(X,\beta)^\T}  \frac{e(\ext^2(\cI,\cI))}{e(\ext^1(\cI,\cI))}.$$
A calculation of the virtual representation
$$\ext^1(\cI,\cI)-\ext^2(\cI,\cI)$$
is required for the
evaluation of the virtual localization formula.

\subsection{Virtual tangent space}

The virtual tangent space at $\mathcal{I}$
is given by
$$
\mathcal{T}_{\left[\mathcal{I}\right]}=
\ext^{1}(\mathcal{I},\mathcal{I}) - \ext^{2}(\mathcal{I},\mathcal{I})
= \chi(\mathcal{O}, \mathcal{O}) - \chi(\mathcal{I}, \mathcal{I})
$$
where
$$
\chi(\mathcal{F}, \mathcal{G}) =
\sum_{i=0}^{3}(-1)^{i}\ext^{i}(\mathcal{F},\mathcal{G})\,.
$$
We can compute each
Euler characteristic using the local-to-global spectral sequence
\begin{align*}
  \chi(\mathcal{I}, \mathcal{I})
&= \sum_{i,j=0}^{3}(-1)^{i+j}H^{i}(\sheafext^{j}(\mathcal{I},\mathcal{I})) \\
&= \sum_{i,j = 0}^{3}(-1)^{i+j}\mathfrak{C}^{i}(\sheafext^{j}(\mathcal{I},\mathcal{I}))
\,,
\end{align*}
where, in the second line,
we have replaced the cohomology terms with the
Cech complex with respect to the open affine cover $\{U_\alpha\}$.
Though these modules are infinite-dimensional,
they have finite-dimensional weight spaces and, therefore,
their $\T$-character is well defined as a formal power
series.

Since $Y$ is supported on the curves $C_{\alpha\beta}$,
we have $\cI = \cO_X$
on the intersection of three or more $U_{\alpha}$.
Therefore, only the $\mathfrak{C}^{0}$ and $\mathfrak{C}^{1}$ terms contribute to
the calculation. We find,
\begin{multline}\label{virtT}
  \mathcal{T}_{\left[\mathcal{I}\right]} =
\bigoplus_{\alpha}\left(\Gamma(U_{\alpha})
- \sum_{i}(-1)^{i}\Gamma(U_{\alpha}, \sheafext^{i}(\mathcal{I},\mathcal{I}))\right) \\
-
\bigoplus_{\alpha, \beta}\left(\Gamma(U_{\alpha\beta}) -
\sum_{i}(-1)^{i}\Gamma(U_{\alpha\beta}, \sheafext^{i}(\mathcal{I},\mathcal{I}))\right)
\,.
\end{multline}
The calculation is reduced to a sum over all the vertices and edges of the
Newton polyhedron.
In each case, we are given an ideal
$$I=I_\alpha,I_{\alpha\beta}\subset \Gamma(U),$$
and we need to compute
$$
\left(\Gamma(U) - \sum_{i}(-1)^{i}Ext^i(I,I)\right)
$$
over the ring $\Gamma(U)$, which is isomorphic to
 $\mathbf{C}[x,y,z]$ in the vertex case and is isomorphic to
$\mathbf{C}[x,y,z,z^{-1}]$ in the
edge case.  We treat each case separately.

\subsection{Vertex calculation}

Let $R$ be the coordinate ring,
$$
R = \C[x_1,x_2,x_3] \cong \Gamma(U_{\alpha}).
$$
As before, we can assume the $\T$-action on $R$ is the
standard action \eqref{st_act}. Consider a $\T$-equivariant graded free
resolution of $I_\alpha$,
\begin{equation}
  \label{resol}
  0 \rightarrow F_{s} \rightarrow \dots \rightarrow F_{2} \rightarrow F_{1} \rightarrow
I_\alpha \rightarrow 0\,,
\end{equation}
such as, for example, the Taylor resolution \cite{Taylor}. Each term in \eqref{resol}
has the form
$$
F_i = \bigoplus_j R(d_{ij})\,, \quad d_{ij} \in \Z^3 \,.
$$
The Poincare polynomial
$$
P_\alpha(t_1,t_2,t_3) = \sum_{i,j} (-1)^i \, t^{d_{ij}}
$$
does not depend on the choice of the resolution
\eqref{resol}. In fact, from
the resolution \eqref{resol} we see
that the Poincare polynomial
$P_\alpha$ is related to the $\T$-character
of $R/I_\alpha$ as follows:
\begin{align}
Q_\alpha(t_1,t_2,t_3)
:\!&= \tr_{R/I_\alpha} (t_1,t_2,t_3) \notag \\
               & = \sum_{(k_1,k_2,k_3)\in \pi_\alpha} t_1^{k_1} t_2^{k_2} t_3^{k_3}
\notag \\
               & = \frac{1+P_\alpha(t_1,t_2,t_3)}{(1-t_1)(1-t_2)(1-t_3)} \label{PQ}\,,
\end{align}
where trace in the first line denotes the trace
of the $\T$-action on $R/I_\alpha$.

The virtual representation
$\chi(I_\alpha,I_\alpha)$ is given by the
following alternating sum
\begin{align*}
\chi(I_\alpha,I_\alpha)  &= \sum_{i,j,k,l} (-1)^{i+k}
\Hom_R(R(d_{ij}), R(d_{kl}))
\\
&=  \sum_{i,j,k,l} (-1)^{i+k}
R(d_{kl}-d_{ij})\,,
\end{align*}
and, therefore,
$$
\tr_{\chi(I_\alpha,I_\alpha)} (t_1,t_2,t_3) =
\frac{P_\alpha(t_1,t_2,t_3) \, P_\alpha(t_1^{-1},t_2^{-1},t_3^{-1})}
{(1-t_1)(1-t_2)(1-t_3)} \,.
$$
We find the character of the $\T$-action
on the $\alpha$ summand of \eqref{virtT} is given by:
$$
\frac{1-P_\alpha(t_1,t_2,t_3) \, P_\alpha(t_1^{-1},t_2^{-1},t_3^{-1})}
{(1-t_1)(1-t_2)(1-t_3)} \,.
$$
Using \eqref{PQ}, we may express the answer in terms of
the generating function $Q_\alpha$ of the partition $\pi_\alpha$,
\begin{multline}\label{vertexchar}
  \tr_{R-\chi(I_\alpha,I_\alpha)}(t_1,t_2,t_3) \\
= Q_{\alpha} -
\frac{\overline{Q}_\alpha}{t_1 t_2 t_3} +  Q_{\alpha}
\overline{Q}_\alpha \frac{(1-t_1)(1-t_2)(1-t_3)}{t_1 t_2 t_3} \,,
\end{multline}
where
$$
\overline{Q}_\alpha(t_1,t_2,t_3) = Q_\alpha(t_1^{-1},t_2^{-1},t_3^{-1}) \,.
$$
The rational function \eqref{vertexchar} should be expanded
in ascending powers of the $t_i$'s.

\subsection{Edge calculation}

We now consider the summand of \eqref{virtT} corresponding to
a pair $(\alpha,\beta)$. Our calculations will involve modules over the ring
$$
R=\Gamma(U_{\alpha\beta})  =
 \C[x_2,x_3]\otimes_\C
\C[x_1,x_1^{-1}] \,.
$$
The $\C[x_1,x_1^{-1}]$ factor will result
only in the overall factor
$$
\delta(t_1) = \sum_{k\in \Z} t_1^k,
$$
the formal $\delta$-function at $t_1=1$, in the $\T$-character.
Let
$$
Q_{\alpha\beta} (t_2,t_3) = \sum_{(k_2,k_3) \in \lambda_{\alpha\beta}} t_2^{k_2}
t_3^{k_3}
$$
be the generating function for the edge partition $\lambda_{\alpha\beta}$.
Arguing as in the vertex case, we find
\begin{multline}\label{edgechar}
 - \tr_{R-\chi(I_{\alpha\beta},I_{\alpha\beta})}(t_1,t_2,t_3) \\
= \delta(t_1) \left(- Q_{\alpha\beta} -
\frac{\overline{Q}_{\alpha\beta}}{t_2 t_3} +  Q_{\alpha\beta}
\overline{Q}_{\alpha\beta} \frac{(1-t_2)(1-t_3)}{t_2 t_3} \right)\,.
\end{multline}
Note that because of the relations
$$
\delta(1/t) = \delta(t) = t\delta(t) \,,
$$
the character \eqref{edgechar} is invariant under the change
of variables \eqref{trans}.

\subsection{The equivariant vertex}
\label{seqv}

The formulas \eqref{vertexchar} and \eqref{edgechar} express
the Laurent polynomial $\tr_{\cT_{\left[\cI\right]}}(t_1,t_2,t_3)$
as a linear combination of infinite formal power series. Our goal
now is to redistribute the terms in these series so that
both the vertex and edge contributions are finite.

The edge character \eqref{edgechar} can be written as
\begin{equation}
  \label{FF}
  \frac{F_{\alpha\beta}(t_2,t_3)}{1-t_1} +
t_1^{-1}\frac{F_{\alpha\beta}(t_2,t_3)}{1-t_1^{-1}} \,,
\end{equation}
where
$$
F_{\alpha\beta}(t_2,t_3) = - Q_{\alpha\beta} -
\frac{\overline{Q}_{\alpha\beta}}{t_2 t_3} +  Q_{\alpha\beta}
\overline{Q}_{\alpha\beta} \frac{(1-t_2)(1-t_3)}{t_2 t_3} \,.
$$
and the first (resp.\ second) term in \eqref{FF} is expanded
in ascending (resp.\ descending) powers of $t_1$.

Let us denote the character \eqref{vertexchar} by $F_\alpha$ and
define
$$
\bV_\alpha = F_\alpha + \sum_{i=1}^3
\frac{F_{\alpha\beta_i}(t_{i'},t_{i''})}{1-t_i}\,,
$$
where $C_{\alpha\beta_1}, C_{\alpha\beta_2}, C_{\alpha\beta_3}$  are the three $\T$-invariant
rational curves passing through the point $X_\alpha \in X^\T$,
and $\{t_i, t_{i'}, t_{i''}\} = \{t_1,t_2,t_3\}$.

Similarly, we define
 \begin{equation*}
   \bE_{\alpha\beta} =
t_1^{-1}\frac{F_{\alpha\beta}(t_2,t_3)}{1-t_1^{-1}}
- \frac{F_{\alpha\beta}\left(t_2\, t_1^{-m_{\alpha\beta}},
t_3\, t_1^{-m'_{\alpha\beta}}\right)}{1-t_1^{-1}}\,.
\end{equation*}
The term $\bE_{\alpha\beta}$ is canonically associated to the edge.
Formulas \eqref{vertexchar} and \eqref{edgechar} yield the following
result.

\begin{Theorem} The $\T$-character of $\cT_{\left[\cI\right]}$
is given by
\begin{equation}
  \label{charTvirt}
  \tr_{\cT_{\left[\cI\right]}}(t_1,t_2,t_3) =
\sum_\alpha \bV_\alpha + \sum_{\alpha\beta} \bE_{\alpha\beta} \,.
\end{equation}
\end{Theorem}

\begin{lm} Both $\bV_\alpha$ and $ \bE_{\alpha\beta}$
are Laurent polynomials.
\end{lm}

\begin{proof}
  The numerator of $\bE_{\alpha\beta}$ vanishes at $t_1=1$,
whence it is divisible by the denominator. The claim for $\bV_\alpha$
follows from
$$
Q_\alpha = \frac{Q_{\alpha\beta}}{1-t_1} + \dots\,,
$$
where the dots stand for terms regular at $t_1=1$.
\end{proof}

{}From $\bV_\alpha$, the equivariant localization formula defines a natural
$3$-parametric family of measures $\bw$ on 3-dimensional
partitions $\pi_\alpha$. Namely, the measure of $\pi_\alpha$ equals
$$
\bw(\pi_\alpha) = \prod_{k\in \Z^3} \left( s,k \right)^{-\bv_k} \,,
$$
where $s=(s_1,s_2,s_3)$ are parameters, $(\, \cdot\, ,\, \cdot\,)$
denotes the standard inner product, and $\bv_k$ is the coefficient
of $t^k$ in $\bV_\alpha$. We call the measure $\bw$
the \emph{equivariant vertex measure}.

\subsection{Local CY and the topological vertex}
We now specialize to the local Calabi-Yau
geometry discussed in Section \ref{ldef}.

Let $S$ be a nonsingular, toric, Fano surface with canonical bundle
$K_S$. We view the total space of $K_S$ as an open toric Calabi-Yau
3-fold.
Let $X$ be the toric compactification defined in Section \ref{ldef}.
By definition,
\begin{equation}
\label{jjj}
\bZ'_{DT}(S;q)_\beta = \bZ_{DT}(X;q)_\beta \big/ \bZ_{DT}(X;q)_0
\,,
\end{equation}
for $\beta \in H_2(S,\Z)$.

We may compute the right side of (\ref{jjj}) by localization.
Let
$$
D= X \setminus S
$$
denote the divisor at infinity.
Let $[\cI]\in I_n(X,\beta)$ be a $\T$-fixed ideal sheaf. We have seen
the weights of the virtual tangent representation of $[\cI]$
are determined by the vertices and edges of the support of $Y$.
Since $\beta$ is a class on $S$, the support of $Y$ lies in $K_S$
except for possibly a finite union of zero dimensional subschemes supported on $D$.
Therefore, as a consequence of the virtual localization formula for
the Donaldson-Thomas theory of $X$,
we find,
\begin{equation}
\label{nnn}
\bZ'_{DT}(S;q)_\beta =
\frac{ \sum_n q^n\sum_{[\cI] \in I_n(K_S,\beta)}  \frac{e(\ext^2(\cI,\cI))}{e(\ext^1(\cI,\cI))}}
     { \sum_n q^n \sum_{[\cI] \in I_n(K_S,0)}  \frac{e(\ext^2(\cI,\cI))}{e(\ext^1(\cI,\cI))}}.
\end{equation}
Here, only the ideal sheaves $\cI$ for which $Y$ has
compact support in $K_S$ are considered.
In particular, the local Donaldson-Thomas theory should
be viewed as independent of the compactification $X$.

The open set $K_S$ has a canonical Calabi-Yau 3-form $\Omega$.
There is 2-dimensional subtorus,
$$ \T_0 \subset \T,$$
which preserves $\Omega$. We will evaluate the formula (\ref{nnn})
on the subtorus $\T_0$.

Let $U_\alpha \subset K_S$ be a chart with coordinates \eqref{st_act}.
The subgroup $\T_0$ is defined by
$$
t_1 t_2 t_3 = 1 \,.
$$
By Serre duality for a compact Calabi-Yau 3-fold, we obtain
a canonical isomorphism
$$
\ext^{1}(\mathcal{I},\mathcal{I})_{0} = \ext^{2}(\mathcal{I},\mathcal{I})_{0}^*.
$$
We will find the $\T_0$-representations to be dual in the local
Calabi-Yau geometry as well. Formula (\ref{nnn}) will be evaluated by
canceling the dual weights and counting signs.

The following functional equation
for the character \eqref{charTvirt} expresses Serre duality. On the
subtorus $t_1 t_2 t_3 =1$, the character is odd
under the involution $f\mapsto \overline f$
defined by
$$
(t_1,t_2,t_3) \mapsto (t_1^{-1},t_2^{-1},t_3^{-1}) \,.
$$
Below we will see, in fact, each term in \eqref{charTvirt}
is an anti-invariant of this transformation.

A crucial
technical point is that no term of \eqref{charTvirt}
specializes to 0 weight under the restriction to $\T_0$.
Since the specializations are all nonzero, the localization
formula for $\T$ may be computed after to restriction to
$\T_0$. We leave the straightforward verification
to the reader.

We will split the edge contributions of  \eqref{charTvirt} in two pieces
$$
\bE_{\alpha \beta}= \bE_{\alpha \beta}^+ +  \bE_{\alpha \beta}^-
$$
satisfying
\begin{equation}
  \label{Feq}
  \overline{\bE}_{\alpha\beta}^+\Big|_{t_1 t_2 t_3 =1}  = -
\bE_{\alpha\beta}^-\Big|_{t_1 t_2 t_3 =1}\,.
\end{equation}
The total count of $(-1)$'s that $\bE_{\alpha \beta}$
contributes is then determined by
the parity of the evaluation of $\bE_{\alpha\beta}^+$
at the point
$(t_1,t_2,t_3)=(1,1,1)$.
Concretely, we set
$$
F_{\alpha\beta}^+ = - Q_{\alpha\beta} -
Q_{\alpha\beta} \overline{Q}_{\alpha\beta}
\frac{1-t_2}{t_2}
$$
and  define $\bE_{\alpha\beta}^+$  in terms of $F_{\alpha\beta}^+$ using
the same formulas as before. A straightforward check
verifies \eqref{Feq}.

Observe that
$$
\bE_{\alpha\beta}^+\Big|_{t_1=1} =
\left(m_{\alpha\beta} \, t_2 \frac{\partial}{\partial t_2} +
m'_{\alpha\beta}\, t_3 \frac{\partial}{\partial t_3} - 1\right) F_{\alpha\beta}^+ \,,
$$
{}Hence, we conclude
\begin{equation}
  \label{parbE}
 \bE_{\alpha\beta}(1,1,1) \equiv f(\alpha,\beta) +
m_{\alpha\beta}|\lambda_{\alpha\beta}| \mod 2 \,,
 \end{equation}
where the function $f(\alpha,\beta)$ was defined in \eqref{fab}. The
second term in \eqref{parbE} comes from applying $\frac{\partial}{\partial t_2}$
to the $(1-t_2)$ factor in the $Q\overline{Q}$-term.

Naively, a similar splitting of the vertex term is given
by defining $F_{\alpha}^+$ to be equal to
$$
Q_\alpha - Q_\alpha \overline{Q}_\alpha
\frac{(1-t_1)(1-t_2)}{t_1 t_2} \,.
$$
However, the definition is not satisfactory since it leads to
rational functions and not polynomials:
the
\begin{equation}
  \label{term}
  \frac{(1-t_1)(1-t_2)(1-t_3)}{t_1 t_2 t_3}
\end{equation}
factor in the $Q\overline{Q}$-term in \eqref{vertexchar}
can be split in
three different ways and no single choice can serve
all terms in the $Q\overline{Q}$-product. The correct
choice of the splitting is the following. Define
the polynomial $Q'_{\alpha}$ by the following
equality
$$
Q_\alpha = Q'_{\alpha} + \sum_{i=1}^3 \frac{Q_{\alpha\beta_i}}{1-t_i}\,.
$$
Now for each set of bar-conjugate terms in the expansion of
the $Q\overline{Q}$-product, we pick its own splitting of
\eqref{term}, so that, for example, the term
$$
\frac{Q_{\alpha\beta} \overline{Q}_{\alpha\beta}}
{(1-t_1)(1-t_1^{-1})}\frac{(1-t_1)(1-t_2)}{t_1 t_2}
$$
cancels the corresponding contribution of $F_{\alpha\beta}^+$, and
for $i\ne j$ the terms
$$
\left( \frac{Q_{\alpha\beta_i} \overline{Q}_{\alpha\beta_j}}
{(1-t_i)(1-t_j^{-1})} +
\frac{\overline{Q}_{\alpha\beta_i} Q_{\alpha\beta_j}}
{(1-t_i^{-1})(1-t_j)}\right)
\frac{(1-t_i)(1-t_j)}{t_i t_j}
$$
are regular and even at $(1,1,1)$.

Using splitting defined above, we easily compute
\begin{equation}
  \label{parbV}
  \bV_\alpha^+(1,1,1) \equiv Q'_{\alpha}(1,1,1) \mod 2 \,.
\end{equation}
{}From the discussion in Section \ref{sdegch}, we find
\begin{equation}
\label{ggg}
Q'_{\alpha}(1,1,1) = |\pi_\alpha| \,.
\end{equation}
Equations \eqref{parbE} and \eqref{ggg} together with Lemma \ref{fchi} yield
the following result.

\begin{Theorem} Let $\cI$ be a $\T$-fixed ideal sheaf in  $I_n(K_S,\beta)$,
$$\frac{e(\ext^2(\cI,\cI))}{e(\ext^1(\cI,\cI))} =
(-1)^{\chi(\oh_Y)+ \sum_{\alpha\beta} m_{\alpha\beta} |\lambda_{\alpha\beta}|}\,,
$$
where the sum in the exponent is over all edges and
$$O(m_{\alpha\beta})\oplus O(m'_{\alpha\beta})$$
is the normal bundle to the edge curve $C_{\alpha\beta}$.
\end{Theorem}

The number $1+m_{\alpha\beta}$ has the same parity as
the \emph{framing} of corresponding edge, which is
a notion introduced
in the context of the topological vertex expansion \cite{topver}.
Comparing the equivariant localization formula for
Donaldson-Thomas theory with the
topological vertex expansion, specifically with
the melting crystal interpretation of the topological
vertex expansion found in \cite{harvard,ORV}, we obtain our main
result.

\begin{Theorem}\label{last}
For the toric local Calabi-Yau geometry, the equality,
$$
\bZ'_{GW}(S;u,v)= \bZ'_{DT}(S;-e^{iu},v)
$$
holds.
\end{Theorem}

The above result depends upon the
evaluation of the Gromov-Witten theory
of  toric local Calabi-Yau 3-folds  via the topological vertex.
The topological vertex has been established in the 1-leg case
in \cite{op,llz} and in the 2-leg case in \cite{llz2}.
These evaluations are sufficient for the GW/DT correspondence for
nonsingular toric local curves and surfaces.


\vspace{+10 pt}
\noindent
Department of Mathematics \\
Princeton University \\
Princeton, NJ 08544, USA\\
dmaulik@math.princeton.edu \\

\vspace{+10 pt}
\noindent
Institut des Hautes Etudes Scientifiques \\
Bures-sur-Yvette, F-91440, France\\
nikita@ihes.fr \\

\vspace{+10 pt}
\noindent
Department of Mathematics \\
Princeton University \\
Princeton, NJ 08544, USA\\
okounkov@math.princeton.edu \\

\vspace{+10 pt}
\noindent
Department of Mathematics\\
Princeton University\\
Princeton, NJ 08544, USA\\
rahulp@math.princeton.edu

\end{document}